\documentclass[12pt]{amsart}
\usepackage{amssymb, amsmath, amscd}

\newtheorem{thm}{Theorem}[section]
\newtheorem{cor}[thm]{Corollary}
\newtheorem{lem}[thm]{Lemma}

\newtheorem{defn}[thm]{Definition}
\newtheorem{rem}[thm]{Remark}

\numberwithin{equation}{section}
\newcommand{\ci}{\circ_{n}}
\newcommand{\be}{\beta}
\newcommand{\ga}{\gamma}
\newcommand{\LL}{L^{\mathcal{E}}}
\newcommand{\LS}{L^{\Se}}
\newcommand{\A}{\mathcal A}

\newcommand{\hh}{\mathfrak{h}}
\newcommand{\C}{\mathbb C}

\newcommand{\W}{\mathcal W}
\newcommand{\One}{{\bf 1}}
\newcommand{\ZZ}{\mathbb Z}
\newcommand{\D}{\mathcal D}
\newcommand{\Se}{\mathcal S}
\newcommand{\cE}{\mathcal E}
\newcommand{\U}{\mathfrak{U}}
\newcommand{\G}{\mathfrak{g}}
\newcommand{\jj}{\mathfrak{j}}
\newcommand{\fa}{\mathfrak{a}}
\newcommand{\bb}{\mathfrak{b}}
\author{Andrew R. Linshaw}
\address{Department of Mathematics\\ Brandeis University\\
Waltham, MA 02454.}
\thanks{I would like to thank Bong H. Lian for many helpful
conversations we have had during the course of this work.}
\title{The Cohomology Algebra
of the Semi-Infinite Weil Complex}
\begin{document}
\maketitle
\begin{abstract}
In 1993, Lian-Zuckerman constructed two cohomology operations on the
BRST complex of a conformal vertex algebra with central charge 26.
They gave explicit generators and relations for the cohomology
algebra equipped with these operations in the case of the $c=1$
model. In this paper, we describe another such example, namely, the
semi-infinite Weil complex of the Virasoro algebra. The
semi-infinite Weil complex of a tame $\mathbb Z$-graded Lie algebra
was defined in 1991 by Feigin-Frenkel, and they computed the linear
structure of its cohomology in the case of the Virasoro algebra. We
build on this result by giving an explicit generator for each
non-zero cohomology class, and describing all algebraic relations in
the sense of Lian-Zuckerman, among these generators.
\end{abstract}

\section{Introduction}
The BRST cohomology of a conformal vertex algebra of central charge
26 is a special case of the semi-infinite cohomology of a tame
$\mathbb Z$-graded Lie algebra $\G$ (in this case the Virasoro
algebra) with coefficients in a $\G$-module $M$. The theory of
semi-infinite cohomology was developed by Feigin and
Frenkel-Garland-Zuckerman \cite{F}\cite{FGZ}, and is an analogue of classical Lie algebra cohomology. In general, there is an
obstruction to the semi-infinite differential being square-zero
which arises as a certain cohomology class in $H^2(\G,\mathbb C)$.
The semi-infinite Weil complex of $\G$ is obtained by taking $M$ to
be the module of \lq\lq adjoint semi-infinite symmetric powers" of
$\G$ \cite{FF1}\cite{AK1}. In this case, an anomaly cancellation
ensures that the differential is always square-zero. The
semi-infinite Weil complex is a vertex algebra and its differential
arises as the zeroth Fourier mode of a vertex operator, a fact which
is useful for doing computations.

This paper is organized as follows. First, we define vertex algebras
and their modules, which have been discussed from various different
points of view in the literature
\cite{B}\cite{FHL}\cite{FLM}\cite{FMS}\cite{K2}\cite{Li}\cite{LZ2}\cite{MS}.
We will follow the formalism introduced in \cite{LZ2}. We describe
the main examples we need, and then define the BRST complex of a
conformal vertex algebra $\A$ with central charge 26. We then recall
the two cohomology operations introduced in \cite{LZ1} on the BRST
cohomology $H^*(\A)$, namely, the {\it dot product} and the {\it
bracket}. We examine in detail the case where $\A = \Se$, ie, the
$\be\ga$-ghost system associated to a one-dimensional vector space.
This coincides with the module of adjoint semi-infinite symmetric
powers of the Virasoro algebra, so the BRST complex of $\Se$ is
exactly the semi-infinite Weil complex of the Virasoro algebra.
Finally, we prove our main result, which is a complete description
of the algebraic structure of $H^*(\Se)$ in the sense of
Lian-Zuckerman.

\begin{thm}
Let $Vir_+$ denote the Lie subalgebra of the Virasoro algebra
generated by $L_n,\ n\ge 0$. As a Lie superalgebra with respect to
the bracket, $H^*(\Se)$ is isomorphic to the semi-direct product of
$Vir_+$ with its adjoint module. As an associative algebra with
respect to the dot product, $H^*(\Se)$ is a polynomial algebra on
one even variable and one odd variable.
\end{thm}

\section{Vertex Algebras}
Let $V=V_0\oplus V_1$ be a super vector space over $\mathbb{C}$, and
let $z,w$ be formal variables. By $QO(V)$, we mean the space of all
linear maps
$$V\rightarrow V((z))=\{\sum_{n\in\mathbb{Z}} v(n) z^{-n-1}|
v(n)\in V,\ v(n)=0\ for\ n>>0 \}.$$ Each element $a\in QO(V)$ can be
uniquely represented as a power series
$a(z)=\sum_{n\in\mathbb{Z}}a(n)z^{-n-1}\in (End\ V)[[z,z^{-1}]]$,
although the latter space is clearly much larger than $QO(V)$. We
refer to $a(n)$ as the $n$-th Fourier mode of $a(z)$. Each $a\in
QO(V)$ is assumed to be of the shape $a=a_0+a_1$ where
$a_i:V_j\rightarrow V_{i+j}((z))$ for $i,j\in\ZZ/2$, and we write
$|a_i| = i$.

On $QO(V)$ there is a set of non-associative bilinear operations,
$\circ_n$, indexed by $n\in\ZZ$, which we call the $n$-th circle
products. They are defined by
\begin{equation}
a(w)\circ_n b(w)=Res_z a(z)b(w)~\iota_{|z|>|w|}(z-w)^n-
(-1)^{|a||b|}Res_z b(w)a(z)~\iota_{|w|>|z|}(z-w)^n.
\end{equation}
Here $\iota_{|z|>|w|}f(z,w)\in\C[[z,z^{-1},w,w^{-1}]]$ denotes the
power series expansion of a rational function $f$ in the region
$|z|>|w|$. Note that $\iota_{|z|>|w|}(z-w)^n \neq \iota_{|w|>|z|}(z-w)^n$ for $n<0$. We usually omit the symbol $\iota_{|z|>|w|}$ and just
write $(z-w)^{n}$ to mean the expansion in the region $|z|>|w|$,
and write $(-1)^n(w-z)^{n}$ to mean the expansion in $|w|>|z|$. It is
easy to check that $a(w)\circ_n b(w)$ above is a well-defined
element of $QO(V)$.

The non-negative circle products are connected through the {\it
operator product expansion} (OPE) formula (\cite{LZ2}, Prop. 2.3).
For $a,b\in QO(V)$, we have \begin{equation} a(z)b(w)=\sum_{n\ge 0}a(w)\circ_n
b(w)~(z-w)^{-n-1}+:a(z)b(w): \end{equation} as formal power series in $z,w$. Here
%(2.2) is customarily written as
%$$a(z)b(w)\sim\sum_{n\ge 0}a(w)\circ_n b(w)~(z-w)^{-n-1},$$ where
%$\sim$ means equal modulo the term $:a(z)b(w):\ $. Here 
$$:a(z)b(w):\ =a(z)_-b(w)\ +\ (-1)^{|a||b|} b(w)a(z)_+\ ,$$ where
$a(z)_-=\sum_{n<0}a(n)z^{-n-1}$ and $a(z)_+=\sum_{n\ge
0}a(n)z^{-n-1}$. (2.2) is customarily written as
$$a(z)b(w)\sim\sum_{n\ge 0}a(w)\circ_n b(w)~(z-w)^{-n-1},$$ where
$\sim$ means equal modulo the term $:a(z)b(w):\ $.

Note that $:a(z)b(z):$ is a well-defined element of
$QO(V)$. It is called the {\it Wick product} of $a$ and $b$, and it
coincides with $a(z)\circ_{-1}b(z)$. The other negative circle products
are related to this by
$$ n!~a(z)\circ_{-n-1}b(z)=\ :(\partial^n a(z))b(z):\ ,$$
where $\partial$ denotes the formal differentiation operator
$\frac{d}{dz}$. For $a_1(z),...,a_k(z)\in QO(V)$, the $k$-fold
iterated Wick product is defined to be
$$:a_1(z)a_2(z)\cdots a_k(z):\ =\ :a_1(z)b(z):$$
where $b(z)=\ :a_2(z)\cdots a_k(z):\ $.

From the definition, we see that \begin{equation}a(z)\circ_0
b(z)=[a(0),b(z)],\end{equation} where the bracket denotes the graded
commutator. It follows that $a\circ_0$ is a graded derivation of
every circle product \cite{LZ2}.

The set $QO(V)$ is a nonassociative algebra with the operations
$\circ_n$ and a unit $1$. We have $1\circ_n a=\delta_{n,-1}a$ for
all $n$, and $a\circ_n 1=\delta_{n,-1}a$ for $n\ge -1$. We are
interested in subalgebras $\A\subset QO(V)$, that is, linear subspaces
of $QO(V)$ containing 1, which are closed under the circle products.
In particular $\A$ is closed under formal differentiation $\partial$
since $\partial a=a\circ_{-2}1$. Following \cite{LZ2}, we call such
a subalgebra a {\it quantum operator algebra} (QOA). Many formal algebraic
notions are immediately clear: a QOA homomorphism is just a linear
map which sends $1$ to $1$ and preserves all circle products; a module over $\A$ is a
vector space $M$ equipped with a QOA homomorphism $\A\rightarrow
QO(M)$, etc. A subset $S=\{a_i|\ i\in I\}$ of $\A$ is said to
generate $\A$ if any element $a\in\A$ can be written as a linear
combination of nonassociative words in the letters $a_i$, $\ci$, for
$i\in I$ and $n\in\mathbb Z$.

\begin{rem}
Fix a nonzero vector $\One\in V$ and let $a,b\in QO(V)$ such that
$a(z)_+\One=b(z)_+\One=0$ for $n\ge 0$. Then it follows immediately
from the definition of the circle products that $(a\circ_p
b)_+(z)\One=0$ for all $p$. Thus if a QOA $\A$ is generated by
elements $a(z)$ with the property that $a(z)_+\One=0$, then every
element in $\A$ has this property. In this case the vector $\One$
determines a linear map
$$
\chi:\A\rightarrow V,~~~a\mapsto a(-1)\One=\lim_{z\rightarrow
0}a(z)\One$$ (called the creation map in \cite{LZ2}), having the
following basic properties: \begin{equation}
\chi(1)=\One,~~~\chi(a\circ_n b)=a(n)b(-1)\One,~~~\chi(\partial^p
a)=p! ~a(-p-1)\One.
\end{equation}
\end{rem}
Next, we define the notion of commutativity in a QOA.
\begin{defn}We say that $a,b\in QO(V)$ circle commute if $$(z-w)^N
[a(z),b(w)]=0$$ for some $N\ge 0$. If $N$ can be chosen to be 0, then
we say that $a,b$ commute. A QOA is said to be commutative if its
elements pairwise circle commute.\end{defn}

The notion of a commutative QOA is abstractly equivalent to the
notion of a vertex algebra (see for e.g. \cite{FLM}). Briefly, every commutative QOA
$\A$ is itself a faithful $\A$-module, called the {\it left regular
module}. Define
$$\rho:\A\rightarrow QO(\A),\ \ \ \ a\mapsto\hat a,\ \ \ \ \hat
a(\zeta)b=\sum_{n\in\mathbb{Z}} (a\circ_n b)~\zeta^{-n-1}.$$ It can
be shown (see \cite{L1} and \cite{LL}) that $\rho$ is an injective QOA homomorphism, and the
quadruple of structures $(\A,\rho,1,\partial)$ is a vertex algebra
in the sense of \cite{FLM}. Conversely, if $(V,Y,\One,D)$ is a
vertex algebra, the subspace $Y(V)\subset QO(V)$ is a commutative
QOA. {\it We will refer to a commutative QOA simply as a vertex
algebra throughout the rest of this paper}.

\begin{rem}
Let $\A'$ be the vertex algebra generated by $\rho(\A)$ inside
$QO(\A)$. Since $\hat a(n)1 = a(z)\ci 1 = 0$ for all $a\in\A$ and
$n\ge 0$, it follows from Remark 2.1 that for every $\alpha\in \A'$,
we have $\alpha_+ 1 = 0$. The creation map $\chi: \A'\rightarrow \A$
sending $\alpha\mapsto \alpha(-1)1$ is clearly a linear isomorphism
since $\chi\circ\rho = id$. It is often convenient to pass between
$\A$ and its image $\A'$ in $QO(\A)$. For example, we shall often
denote the Fourier mode $\hat{a}(n)$ simply by $a(n)$. When we say
that a vertex operator $b(z)$ is annihilated by the Fourier mode
$a(n)$ of a vertex operator $a(z)$, we mean that $a\ci b = 0$. Here
we are regarding $b$ as an element of the state space $\A$, while
$a$ operates on the state space, and the map $a\mapsto \hat{a}$ is
the state-operator correspondence.
\end{rem}

For later use, we write down a formula, valid in any vertex
algebra, which measures the non-associativity of the Wick product.
\begin{lem} Let $\A$ be a vertex algebra, and let $a,b,c\in \A$. Then
\begin{equation}:(:ab:)c:\  - \ :abc:\
 = \sum_{n\ge 0}\frac{1}{(n+1)!}\bigg(:(\partial^{n+1} a)(b\ci c):\ +
(-1)^{|a||b|} :(\partial^{n+1} b)(a\ci c):\bigg)\ .\end{equation}
Note that this sum is finite by circle commutativity. In
particular, we see that $:(:ab:)c:$ and $:abc:$ differ by terms of
the form $:(\partial^i a)X:$ and $:(\partial^i b)Y:$, where $i\ge
1$ and $X,Y \in \A$.
\end{lem}
\begin{proof} By the preceding remark, it suffices to show that
$\hat{a},\hat{b},\hat{c}$ satisfy this identity, which can be
checked by applying the creation map to both sides and then using
(2.4).
\end{proof}

\subsection{Virasoro elements} Many vertex algebras $\A$ have a
{\it Virasoro element}, that is, a vertex operator
$L(z)=\sum_{n\in\mathbb{Z}} L(n)z^{-n-1}$ satisfying the OPE
\begin{equation} L(z)L(w) \sim
\frac{k}{2}\ (z-w)^{-4} + 2L(w)(z-w)^{-2} + \partial
L(w)(z-w)^{-1},
\end{equation}
where the constant $k$ is called the {\it central charge}. It is
customary to write $L(z) = \sum_{n\in\mathbb{Z}} L_n z^{-n-2}$,
where $L_n := L(n+1)$. The Fourier modes $\{L_n|\
n\in\mathbb{Z}\}$ together with a central element $\kappa$ then
generate a copy of the Virasoro Lie algebra $Vir$:
$$[L_n, L_m] = (n-m)L_{n+m} +
\frac{1}{12}(n^3-n)\delta_{n,-m}\kappa.$$ 

Often we require further that $L_0$ be diagonalizable and $L_{-1}$ acts on $\A$ by formal
differentiation. In this case, $(\A,L)$ is called a {\it conformal
vertex algebra}. An element $a\in \A$ which is an eigenvector of
$L_0$ with eigenvalue $\Delta\in\mathbb{C}$ is said to have
conformal weight $\Delta$. In any conformal vertex algebra, the
operation $\ci$ is homogeneous of conformal weight $-n-1$. In
particular, the Wick product $\circ_{-1}$ is homogeneous of
conformal weight zero.

\subsection{$\beta\gamma$- and $bc$-ghost systems}
Let $V$ be a finite-dimensional vector space. Regard $V\oplus V^*$
as an abelian Lie algebra. Then its loop algebra has a
one-dimensional central extension by $\mathbb{C}\tau$
$$\mathfrak{h} = \hh(V) = (V\oplus V^*)[t,t^{-1}]\oplus \mathbb{C}\tau,$$
which is known as a Heisenberg algebra. Its bracket is given by
$$[(x,x')t^n,(y,y')t^m]=(\langle y',x\rangle-\langle
x',y\rangle)\delta_{n+m,0}\tau,$$ for $x,y\in V$ and $x',y'\in V^*$. Let $\bb\subset\hh$ be the subalgebra generated by $\tau$, $(x,0)t^n$, and $(0,x')t^m$, for $n\ge 0$ and $m>0$, and let $C$
be the one-dimensional $\bb$-module on which $(x,0)t^n$ and
$(0,x')t^m$ act trivially and the central element $\tau$ acts by
the identity. Denote the linear operators representing
$(x,0)t^n,(0,x')t^n$ on $\U\hh\otimes_{\U\bb}C$ by
$\beta^x(n),\gamma^{x'}(n-1)$, respectively, for $n\in\mathbb{Z}$.
The power series
$$\beta^x(z)=\sum_{n\in\mathbb{Z}}\beta^x(n)z^{-n-1},\ \ \ 
\gamma^{x'}(z)=\sum_{n\in\mathbb{Z}}\gamma^{x'}(n)z^{-n-1}\ \in
QO(\U\hh\otimes_{\U\bb}C)$$ generate a vertex algebra $\Se(V)$
inside $QO(\U\hh\otimes_{\U\bb}C)$, and the generators satisfy the OPE relations
$$\beta^x(z)\gamma^{x'}(w)\sim\langle x',x\rangle (z-w)^{-1},\ \ \ \beta^x(z)\beta^y(w)\sim 0,\ \ \ \gamma^{x'}(z)\gamma^{y'}(w)\sim 0.$$
This algebra was introduced in \cite{FMS} and is
known as a $\beta\gamma$-ghost system, or a semi-infinite symmetric algebra. The creation map $$\chi: \Se(V)\rightarrow
\U\hh\otimes_{\U\bb}C,$$ which sends $a(z)\mapsto a(-1)(1\otimes 1)$,
is easily seen to be a linear isomorphism. By the
Poincare-Birkhoff-Witt theorem, the vector space
$\U\hh\otimes_{\U\bb}C$ has the structure of a polynomial algebra
with generators given by the negative Fourier modes
$\beta^x(n),\gamma^{x'}(n)$, $n<0$, which are linear in $x\in V$ and
$x'\in V^*$. It follows from (2.4) that $\Se(V)$ is spanned by the
collection of iterated Wick products of the form
$$\mu=\ :\partial^{n_1}\be^{x_1}\cdots
\partial^{n_s}\be^{x_s}
\partial^{m_1}\ga^{x'_1}\cdots \partial^{m_t}\ga^{x'_t}:\ .$$

$\Se(V)$ has a natural $\mathbb{Z}$-grading which we call the
$\be\ga$-ghost number. Fix a basis $x_1,\dots,x_n$ for $V$ and a
corresponding dual basis $x'_1,\dots,x'_n$ for $V^*$. Define the
$\be\ga$-ghost number to be the eigenvalue of the diagonalizable
operator $[B,-]$, where $B$ is the zeroth Fourier mode of the vertex
operator
$$\sum_{i=1}^n:\be^{x_i}\ga^{x'_i}:\ .$$ Clearly $B$ is
independent of our chosen basis of $V$, and $\beta^x, \gamma^{x'}$ have $\be\ga$-ghost numbers $-1,1$ respectively.

We can also regard $V\oplus V^*$ as an odd abelian Lie (super)
algebra, and consider its loop algebra and a one-dimensional
central extension by $\mathbb{C}\tau$ with bracket
$$
[(x,x')t^n,(y,y')t^m]=(\langle y',x\rangle+\langle
x',y\rangle)\delta_{n+m,0}\tau.
$$

Call this Lie algebra $\jj=\jj(V)$, and form the induced module
$\U\jj\otimes_{\U\fa}C$. Here $\fa$ is the subalgebra of $\jj$
generated by $\tau$, $(x,0)t^n$, and $(0,x')t^m$, for $n\geq 0$ and
$m>0$, and $C$ is the one-dimensional $\fa$-module on which
$(x,0)t^n$ and $(0,x')t^m$ act trivially and $\tau$ acts by $1$.
There is a vertex algebra $\cE(V)$, analogous to $\Se(V)$, which is
generated by the odd vertex operators
$$b^x(z)=\sum_{n\in\mathbb{Z}} b^x(n)z^{-n-1},\ \ \
c^{x'}(z)=\sum_{n\in\mathbb{Z}}c^{x'}(n)z^{-n-1}\ \in
QO(\U\jj\otimes_{\U\fa}C),$$ which satisfy the OPE relations
$$
b^x(z)c^{x'}(w)\sim\langle x',x\rangle(z-w)^{-1},\ \ \ b^x(z)b^y(w)\sim 0,\ \ \ c^{x'}(z)c^{y'}(w)\sim 0.$$ 

This vertex algebra is known as a $bc$-ghost system, or a semi-infinite
exterior algebra. Again the creation map
$\cE(V)\rightarrow \U\jj\otimes_{\U\fa}C$, $a(z)\mapsto
a(-1)(1\otimes 1)$, is a linear isomorphism. As in the symmetric
case, the vector space $\U\jj\otimes_{\U\fa}C$ has the structure
of an odd polynomial algebra with generators given by the negative
Fourier modes $b^x(n),c^{x'}(n)$, $n<0$, which are linear in $x\in
V$ and $x'\in V^*$. As above, it follows that $\cE(V)$ is spanned
by the collection of all iterated Wick products of the vertex
operators $\partial^k b^x$ and $\partial^k c^{x'}$, for $k\ge 0$.

$\cE(V)$ has a $\mathbb{Z}$-grading which we call the $bc$-ghost
number (or fermion number). It is given by the eigenvalue of the
diagonalizable operator $[F,-]$, where $F$ is the zeroth Fourier
mode of the vertex operator $$-\sum_{i=1}^n :b^{x_i}c^{x'_i}:\ .$$ $F$ is independent of our choice of basis for $V$, and $b^x,c^{x'}$ have $bc$-ghost numbers $-1,1$ respectively.
We will denote the $bc$-ghost number of a homogeneous element $a\in\cE(V)$ by $|a|$. Note that this coincides with our earlier notation for the $\mathbb{Z}/2$-grading on $\cE(V)$ coming from its vertex
superalgebra structure. This causes no difficulty; since
$b^x,c^{x'}$ are odd vertex operators, the $mod\ 2$ reduction of the
$bc$-ghost number coincides with this $\mathbb{Z}/2$-grading.

Let us specialize to the case where $V$ is a one-dimensional
vector space. In this case, $\Se(V)$ coincides with the module of
adjoint semi-infinite symmetric powers of the Virasoro algebra
\cite{FF1}. Fix a basis element $x$ of $V$ and a dual basis
element $x'$ of $V^*$. We denote $\Se(V)$ by $\Se$, and we denote
the generators $\be^x,\ga^{x'}$ by $\be,\ga$, respectively.
Similarly, we denote $\cE(V)$ by $\cE$, and we denote the
generators $b^x,c^{x'}$ by $b,c$, respectively.

For a fixed scalar $\lambda\in\mathbb{C}$, define
\begin{equation} L^{\Se}_{\lambda} = (\lambda-1):\partial
\be\ga:\ +\lambda :\be\partial \ga:\ \in \Se.\end{equation}
An OPE calculation shows that
$$L^{\Se}_{\lambda}(z)\be(w)\sim \lambda \be(w)(z-w)^{-2} + \partial \be(w)(z-w)^{-1},$$
$$L^{\Se}_{\lambda}(z)\ga(w) \sim(1-\lambda)\ga(w)(z-w)^{-2} + \partial \ga(w)(z-w)^{-1},$$
$$ L^{\Se}_{\lambda}(z)L^{\Se}_{\lambda}(w)\sim \frac{k}{2}(z-w)^{-4} +
2L^{\Se}_{\lambda}(w)(z-w)^{-2} + \partial
L^{\Se}_{\lambda}(w)(z-w)^{-1},$$
where $k = 12\lambda^2 - 12\lambda + 2$. Hence $L^{\Se}_{\lambda}$
is a Virasoro element of central charge $k$, and
$(\Se,L^{\Se}_{\lambda})$ is a conformal vertex algebra in which
$\be, \ga$ have conformal weights $\lambda, 1-\lambda$
respectively.

Similarly, define \begin{equation} L^{\cE}_{\lambda} = (1-\lambda):\partial
bc: - \lambda :b\partial c:\ \in \cE.\end{equation} A calculation shows that
$L^{\cE}_{\lambda}$ is a Virasoro element with central charge $k =
-12\lambda^2 + 12 \lambda -2$, $(\cE,L^{\cE}_{\lambda})$ is a
conformal vertex algebra, and $b,c$ have conformal weights
$\lambda, 1-\lambda$ respectively.

\section{BRST cohomology} Observe that if $\A, \A'$ are conformal
vertex algebras with Virasoro elements $L^{\A}, L^{\A'}$ of central
charges $k, k'$, respectively, then $\A\otimes \A'$ is a conformal
vertex algebra with Virasoro element $L^{\A\otimes \A'} = L^{\A} +
L^{\A'}$ (ie, $L^{\A} \otimes 1 + 1 \otimes L^{\A'})$ of central
charge $k + k'$. To simplify notation, the ordered product $ab
=a(z)b(z)$ of two vertex operators $a,b$ in the same formal variable
$z$ will always denote the Wick product.

Fix $\lambda = 2$ in (2.8), and denote the corresponding Virasoro element
$L^{\cE}_2 = -\partial bc - 2b\partial c \in \cE$, by $\LL$. With
this choice, $(\cE,\LL)$ is a conformal vertex algebra of central
charge -26. For any conformal vertex algebra $(\A,L^{\A})$ of
central charge $k$, let $C^*(\A) = \cE \otimes \A$. Denote the
Virasoro element $\LL + L^{\A}$, by $L^C$. The conformal weight
and $bc$-ghost number are given, respectively, by the eigenvalues
of the operators $[L^C_0,-]$ and $[F\otimes 1,-]$ on $C^*(\A)$.

\begin{defn} Let $J_{\A}$ be the following element of $C^*(\A)$:
\begin{equation} J_{\A} = (L^{\A} + \frac{1}{2}\LL)c +
\frac{3}{4}\partial^2c
\end{equation}
\end{defn}
A calculation shows that
\begin{equation}
J_{\A}(z)\circ_0 b(z) = L^C(z).
\end{equation}

We will denote the zeroth Fourier mode $J_{\A}(0)$ by $Q$, so we may rewrite this equation as $ [Q,b(z)] = L^C(z)$ by (2.3). Note that the operator $[Q,-]$ preserves
conformal weight and raises $bc$-ghost number by $1$.

\begin{lem} $Q^2 = 0$ iff $k = 26$. In this case, we can consider
$C^*(\A)$ to be a cochain complex graded by $bc$-ghost number,
with differential $[Q,-]$. Its cohomology is called the BRST
cohomology associated to $\A$, and will be denoted by $H^*(\A)$.
\end{lem}
\begin{proof} First, note that $Q^2 = \frac{1}{2}[Q,Q] = \frac{1}{2}Res_w
J_{\A}(w)\circ_0 J_{\A}(w)$. Computing the OPE of
$J_{\A}(z)J_{\A}(w)$ and extracting the coefficient of
$(z-w)^{-1}$, we find that $J_{\A}(w) \circ_0 J_{\A}(w) =
\frac{3}{2}\partial (\partial^2c(w)c(w)) +
\frac{k-26}{12}\partial^3c(w)c(w)$. Since the residue of a total
derivative is zero, only the second term contributes, and it
follows that $Res_w J_{\A}(w)\circ_0 J_{\A}(w) = 0$ iff
$k=26$.\end{proof} From now on, we will only consider the case
where $k=26$.

\section{Algebraic Structure of $H^*(\A)$}
In this section, we recall without proof some facts from
\cite{LZ1} on the algebraic structure of the BRST cohomology. We
first note that any cohomology class can be represented by an
element $u(z)$ of conformal weight $0$, since (3.2) implies that
$[Q,b(1)] = L^C_0$. Since $[Q,-]$ acts by derivation on
each of the products $\ci$ on $C^*(\A)$, each $\ci$ descends to a
product $H^*(\A)$. Since $\ci$ lowers conformal weight by $n+1$,
all these products are trivial except for the one induced by
$\circ_{-1}$ (the Wick product), which we call the {\it dot
product}. We write the dot product of $u$ and $v$ as $uv$.

The cohomology $H^*(\A)$ has another bilinear operation known as the
{\it bracket}. First, we define the bracket on the space $C^*(\A)$.
\begin{defn} Given $u(z),v(z) \in C^*(\A)$, let
\begin{equation}\{u(z),v(z)\} = (-1)^{|u|}(b(z)\circ_0
u(z))\circ_0 v(z).
\end{equation}
\end{defn}
The equivalence between this definition and the one given in
\cite{LZ1} is shown in \cite{LZ2}. From this description, it is
easy to see that the bracket descends to $H^*(\A)$, inducing a
well-defined bilinear operation.
\begin{thm}
The following algebraic identities hold on $H^*(\A)$:
\begin{equation} uv = (-1)^{|u||v|}vu \end{equation}
\begin{equation}(uv)t = u(vt)\end{equation}
\begin{equation}\{u,v\}= -(-1)^{(|u|-1)(|v|-1)}\{v,u\}\end{equation}
\begin{equation}(-1)^{(|u|-1)(|t|-1)}\{u,\{v,t\}\}\ +\
(-1)^{(|t|-1)(|v|-1)}\{t,\{u,v\}\}\ +\
(-1)^{(|v|-1)(|u|-1)}\{v,\{t,u\}\}=0\end{equation}
\begin{equation}\{u,vt\} = \{u,v\}t + (-1)^{(|u|-1)(|v|)}v
\{u,t\}\end{equation}
\begin{equation}b\circ_1 \{u,v\} = \{b\circ_1 u,v\} +
(-1)^{|u|-1}\{u,b\circ_1 v\}\end{equation}
\begin{equation}\{,\}:H^p\times H^q \to H^{p+q-1}\end{equation}
(4.2)-(4.3) say that $H^*(\A)$ is an associative, graded
commutative algebra with respect to the dot product. (4.4)-(4.5)
say that under the bracket, $H^*(\A)$ is a Lie superalgebra with
respect to the grading ($bc$-ghost number - 1). Also, note that
$H^1(\A)$ is an ordinary Lie algebra under the bracket. Taking
$p=1$ in (4.8), we see that for every $q$, $H^q(\A)$ is a module
over $H^1(\A)$. \end{thm}

\section{The Semi-infinite Weil Complex of the Virasoro Algebra}
In this section, we give a complete description of the algebraic
structure of $H^*(\A)$ in the case $\A = \Se$, with the choice
$\lambda = 2$ in (2.7). In this case, the Virasoro element $\LS = L^{\Se}_2
=\partial \be\ga + 2\be\partial \ga \in \Se$ has central charge
26.
\begin{defn}$C^*(\Se) = \cE \otimes \Se$, equipped with the
differential $[Q,-]$ and the Virasoro element $L^{\W} = \LL + \LS$,
is called the semi-infinite Weil complex associated to $Vir$, and
will be denoted by $\W$.
\end{defn}
Note that $\W$ is naturally {\it triply} graded. In addition to the
conformal weight and $bc$-ghost number, $\W$ is graded by the
$\be\ga$-ghost number, which is the eigenvalue of $[1\otimes B,-]$.
Note that $[Q,-]$ preserves the $\be\ga$-ghost number. Let $\W^{i,j}
\subseteq \W$ denote the conformal weight zero subspace of
$bc$-ghost number $i$, $\be\ga$-ghost number $j$, and let $Z^{i,j},
B^{i,j} \subseteq \W^{i,j}$ denote the cocycles and coboundaries,
respectively, with respect to $[Q,-]$. Let $H^{i,j} =
Z^{i,j}/B^{i,j}$. Note that $H^i(\Se)$ decomposes as the direct sum
$\bigoplus_{j\in \mathbb Z}H^{i,j}$.

In \cite{FF1} and \cite{FF2}, Feigin-Frenkel computed the {\it
linear} structure of $H^*(\Se)$, namely, the dimension of each of
the spaces $H^{i,j}$.

\begin{thm} For all $j\ge 0$, $dim\ H^{0,j} = dim\ H^{1,j} = 1$.
For all other values of $i,j$, $dim\ H^{i,j} =0$.
\end{thm}

This was proved by using the Friedan-Martinec-Shenker bosonization
\cite{FMS} to express $\Se$ as a submodule of a direct sum of
Feigin-Fuchs modules over $Vir$, and then using known results on the
structure of these modules. We will assume the results in \cite{FF1}
and \cite{FF2}, and use them to describe the {\it algebraic}
structure of $H^*(\Se)$. Our first step is to find a canonical
generator for each non-zero cohomology class. Recall that $\W$ has
a basis consisting of the monomials:
\begin{equation}
\partial^{n_1}b\cdots \partial^{n_i}b\
\partial^{m_1}c\cdots \partial^{m_j}c\
\partial^{s_1}\be\cdots \partial^{s_k}\be\
\partial^{t_1}\ga\cdots \partial^{t_l}\ga
\end{equation}
with $n_1>...> n_i \ge 0$, $ m_1>...> m_j\ge 0$ and $s_1\ge...\ge
s_k \ge 0$, $ t_1\ge...\ge t_l\ge 0$. Let $\mathcal D \subset \W$
be the subspace spanned by monomials which contain at least one
derivative, ie, at least one of the numbers $n_1,...n_i,
m_1,...m_j, s_1,...,s_k, t_1,...,t_l$ above is positive.

\begin{lem} The image of $[Q,-]$ is contained in $\mathcal D$.
\end{lem}
\begin{proof} A straightforward calculation shows that
$$[Q,b] = L^{\W},\ \ [Q,c] = c\partial c,\ \ [Q,\be] = c\partial \be
+ 2\partial c \be,\ \ [Q,\ga] = c\partial \ga -\partial c \ga.$$ 
The claim follows by applying the graded derivation $[Q,-]$ to a
monomial of the form (5.1), and then using (2.5) to express the
result as a linear combination of standard monomials of the form
(5.1). Note that for any vertex operators $a,b,c\in \W$, the
expression $:(ab:)c:\ -\ :abc:\ $, which measures the
non-associativity of the Wick product, always lies in
$\D$.\end{proof}

\begin{lem} Let $x=\be\ga^2-bc\ga + \frac{3}{2}\partial\ga$. Then
$x\in Z^{0,1}$ and $x \notin B^{0,1}$, so $x$ represents a
non-zero cohomology class. Since $H^{0,1}$ is $1$-dimensional, $x$
generates $H^{0,1}$. Similarly, let $y = c\be\ga +
\frac{3}{2}\partial c$. Then $y\in Z^{1,0}$ and $y\notin B^{1,0}$,
so $y$ generates $H^{1,0}$.
\end{lem}
\begin{proof} The proof that $x\in Z^{0,1}$ and $y\in Z^{1,0}$ is a
straightforward calculation. Since $x$ contains a monomial with no
derivatives, $x\notin \mathcal D$. By Lemma 5.3, $B^{0,1} \subset
\mathcal D$, so $x \notin B^{0,1}$. Similarly, $y\notin
B^{1,0}$.\end{proof}

Our main result is the following
\begin{thm} For each integer $k\ge 0$, $x^k$ represents a non-zero cohomology
class in $H^{0,k}$, and $yx^k$ represents a non-zero class in
$H^{1,k}$. By Theorem 5.2, these are all the non-zero classes in
$H^*(\Se)$.
\end{thm}
\begin{proof}It is clear from the derivation property of $[Q,-]$ that
$x^k\in Z^{0,k}$ and $yx^k\in Z^{1,k}$, so it suffices to show that
$x^k\notin B^{0,k}$ and $yx^k\notin B^{1,k}$. For each integer $k
\ge 0$, define: $$x_k =\be^k\ga^{2k}- kbc\be^{k-1}\ga^{2k-1},\ \ \ \
\ \ \ \ y_k = c\be^{k+1}\ga^{2k+1}.$$ We claim that:
\begin{equation}x^k = x_k + D_k,
\end{equation}
\begin{equation}yx^k = y_k + D'_k, \end{equation}
for some $D_k\in\mathcal{D}$ and $D'_k\in\mathcal{D}$. Since $x_k$
and $y_k$ have no derivatives, it follows that $x^k$ and $yx^k$ do
not lie in $\mathcal{D}$. Now we can apply Lemma 5.3 to conclude
that $x^k\notin B^{0,k}$ and $yx^k\notin B^{1,k}$.

We begin with (5.2) and proceed by induction. The cases $k=0$ and
$k=1$ are obvious, so assume the statement true for $k-1$.
\begin{equation}x^k = (\be\ga^2-bc\ga+\frac{3}{2}\partial
\ga)x^{k-1}\end{equation} $$ =
(\be\ga^2-bc\ga)(\be^{k-1}\ga^{2k-2}-(k-1)bc\be^{k-2}\ga^{2k-3}+
D_{k-1}) +E_0,$$ where $E_0=\frac{3}{2}\partial\ga x^{k-1}.$ We
expand this product and apply Lemma 2.4 repeatedly:
\begin{equation}
(bc\ga)((k-1)bc\be^{k-2}\ga^{2k-3}) =
(k-1)b^2c^2\be^{k-2}\ga^{2k-2} + E_1 = E_1
\end{equation} since $b,c$ are anti-commuting variables.
\begin{equation} (bc\ga)(\be^{k-1}\ga^{2k-2}) = bc\be^{k-1}\ga^{2k-1}
+ E_2
\end{equation}
\begin{equation} (\be\ga^2)((k-1)bc\be^{k-2}\ga^{2k-3}) =
(k-1)bc\be^{k-1}\ga^{2k-1} + E_3
\end{equation}
\begin{equation} (\be\ga^2)(\be^{k-1}\ga^{2k-2}) = \be^k\ga^{2k} + E_4
\end{equation}
where $E_i\in \mathcal D$ for $i = 0,1,2,3,4$. It is easy to see
that $(bc\ga) D_{k-1}\in \D$ and $(\be\ga^2)D_{k-1}\in \mathcal D$.
(5.4) follows by collecting terms from (5.5)-(5.8). Finally, the
same argument proves (5.3). \end{proof}

Using Theorem 5.5, we can now describe the algebraic structure of
$H^*(\Se)$. Let $Vir_+\subseteq Vir$ be the Lie subalgebra generated
by $L_n, n\ge 0$. The Cartan subalgebra $h$ of $Vir_+$ is generated
by $L_0$.
\begin{cor} As an associative algebra with respect to the dot product, $H^*(\Se)$ is a polynomial
algebra in one even variable, $x$, and one odd variable, $y$. In
other words, $H^*(\Se)$ is isomorphic to the classical Weil algebra
associated to $h$.
\end{cor}
It is easy to check from the definition of the bracket (4.1) that
$\{y,x\} = -x$. Using the graded derivation property of the
bracket with respect to the dot product, we can write down all the
bracket relations in $H^*(\Se)$. For any $n,m \ge 0$,
\begin{equation}
\{x^n,x^m\} = 0,\ \ \
 \{yx^n,x^m\} = (n-m)x^{n+m},\ \ \
 \{yx^n,yx^m\} = (n-m)yx^{n+m}.
\end{equation}
It follows that as a Lie algebra, $H^{1,*}$ is isomorphic to $Vir_+$
under the isomorphism $yx^k \mapsto L_k, k\ge 0$. As an
$H^{1,*}$-module, $H^{0,*}$ is isomorphic to the adjoint
representation of $Vir_+$. Finally, we obtain
\begin{cor}
As a Lie superalgebra with respect to the bracket, $H^*(\Se)$ is
isomorphic to the semi-direct product of $Vir_+$ with its adjoint
module.
\end{cor}

\end{document}